\documentstyle[amssymb,amscd]{amsart}
\input xy
\xyoption{all}
\begin{document}
\thispagestyle{empty}
\newcommand{\Aut}{\operatorname{Aut}}

\begin{center}
 {\Large\bf On the AF Embeddability of Crossed Products of AF 
            Algebras by the Integers}
\end{center}
\vspace{5mm}

\begin{center}{\bf Nathanial P. Brown}\\
              {\bf Purdue University}
\end{center}

\vspace{5mm}

\begin{abstract}
 {\small  If $A$ is an AF algebra and $\alpha \in \Aut(A)$, it is
          shown that AF embeddability of the crossed product, 
          $A\times_{\alpha}\Bbb Z$, is equivalent to 
          $A\times_{\alpha}\Bbb Z$ being stably finite.
          This equivalence follows from a simple K-theoretic characterization 
          of AF embeddability. }
\end{abstract}

\vspace{4mm}

     This paper is concerned with the question of when the crossed product of 
an AF algebra by an action of $\Bbb Z$ is itself AF embeddable.  It is well 
known that quasidiagonality and stable finiteness are hereditary properties. 
That is, if $A$ and $B$ are C$^{*}$-algebras with $A \subset B$ and $B$ has 
either of these properties, then so does $A$.  Since AF algebras enjoy both 
of these properties we have that quasidiagonality and stable finiteness are 
geometric obstructions to AF embeddability.  For crossed products of AF 
algebras by $\Bbb Z$, these turn out to be the only obstructions.\\

     If $A$ is an AF algebra, then we may easily describe an algebraic 
obstruction to the AF embeddability of $A\times_{\alpha}\Bbb Z$.\\

{\bf Definition 0.1}  If $A$ is an AF algebra and $\alpha \in \Aut(A)$ then we 
    denote by $H_{\alpha}$ the subgroup of $K_{0} (A)$ given by all 
    elements of the form $\alpha_{*}(x) - x$ for $ x \in K_{0} (A)$. \\

     It follows from the Pimsner-Voiculescu six term exact sequence ([PV]) 
that if $A$ is unital and AF then $K_{0} (A\times_{\alpha}\Bbb Z$) = 
$K_{0} (A)/H_{\alpha}$.  Now, if $B$ is unital and stably finite 
(in particular, if 
$B$ is AF) and $p \in M_{n}(B)$ is a projection then $[p]$ must be a 
nonzero element of $K_{0} (B)$.  Thus, if $A\times_{\alpha}\Bbb Z$ 
embeds into 
B (or if $A\times_{\alpha}\Bbb Z$ is already stably finite) then 
every projection in the matrices over $A\times_{\alpha}\Bbb Z$ 
must give a nonzero element in $K_0 (A\times_{\alpha}\Bbb Z)$.  In 
particular, since $A \hookrightarrow A\times_{\alpha}\Bbb Z$, and we 
have observed that $K_0 (A\times_{\alpha}\Bbb Z) = K_0(A)/H_{\alpha}$, 
we then conclude that $H_{\alpha} \cap K_{0}^{+} (A) = \{0\}$, where 
$K_{0}^{+} (A)$ is the positive cone of $K_{0} (A)$.  Thus, an
algebraic obstruction to AF embeddability of a crossed product would be 
if $H_{\alpha}$ contained a positive element of $K_0 (A)$.  It turns out 
that (so long as $A$ is AF) this is the only algebraic obstruction.  
Indeed, the main result of 
this paper is the following.\\

{\bf Theorem 0.2} If $A$ is an AF algebra (not necessarily unital) and 
$\alpha \in \Aut (A)$ then the following are equivalent.

\begin{enumerate}
  \item $A\times_{\alpha}\Bbb Z$ is AF embeddable
  \item $A\times_{\alpha}\Bbb Z$ is quasidiagonal
  \item $A\times_{\alpha}\Bbb Z$ is stably finite
  \item $H_{\alpha} \cap K_{0}^{+} (A) = \{0\}$
\end{enumerate}

     The implications $(1) \Rightarrow (2) \Rightarrow (3) \Rightarrow (4)$ 
are relatively straightforward (whether or not $A$ is unital).  In fact, 
if $A$ is any C$^{*}$-algebra with a countable approximate unit consisting 
of projections, then the implications $(1) \Rightarrow (2) \Rightarrow 
(3) \Rightarrow (4)$ hold.  However 
the remaining implication, $(4) \Rightarrow (1)$, will 
require some technical machinery (see Section 3 for the details of all 
the implications).  This paper consists of five sections.

     In Section 1 we collect a few useful facts and set some of our notation.

\vspace{1mm}

     In Section 2 we develop the key technical tools.  Using the Rohlin 
property for automorphisms of AF algebras (see Definition 2.1) we will 
show that certain commutative diagrams at the level of K-theory lift to 
commutative diagrams on the algebras.  In recent years, the 
Rohlin property has been studied intensively.  There are notions of the 
Rohlin property for more general group actions ([Oc], [Na]) and this 
notion for actions of $\Bbb Z$  has been used in many 
interesting applications ([R$\o$], [Ki2], [BKRS], [Co2], [EK], just 
to name a few).  All of 
the embedding results presented here could be viewed as more applications 
of the Rohlin property.

\vspace{1mm}

     In Section 3 we prove a  useful K-theoretic characterization of 
AF embeddability (different than that given in Theorem 0.2) and use this 
to prove  Theorem 0.2.

\vspace{1mm}

     In Section 4 we present several applications of our results.  Among 
other things, we will prove a natural generalization of a 
result of Voiculescu that states that if some power of an 
automorphism of an AF algebra is approximately inner then the 
corresponding crossed 
product is AF embeddable (see Theorem 3.6 in [Vo]).
We will also recover Pimsner's criterion for AF embeddability 
of crossed products of $C(X)$ by $\Bbb Z$, in the case $X$ is compact, 
metrizable 
and totally disconnected. 

\vspace{1mm}

     In Section 5 we show that our constructions can be done in such a way 
as to yield rationally injective maps on $K_0(A\times_{\alpha}\Bbb Z)$ and 
thus injective maps when $K_0(A\times_{\alpha}\Bbb Z)$ is assumed to be 
torsion free.

\vspace{1mm}

     The key technical results presented here were established 
 while the author participated in the 
Monbusho Summer Program for Young Foreign Researchers in Sapporo, Japan. 
The author would like to thank Monbusho and the NSF for sponsoring this 
program and allowing him the opportunity to participate in it.  The 
author is particularly indebted to his host, Professor Akitaka Kishimoto 
(nihongo de Kishimoto Sensei), for the many invaluable discussions on 
crossed products, dimension groups and most importantly, the Rohlin 
property. Indeed, these notes owe a great deal to the experience and 
insights shared with the author during his stay in Sapporo.  The author 
would also like to thank Professor Hiroshi Takai for a very profitable
discussion at Tokyo Metropolitan University.  A few of the results presented 
here are (partial) answers to questions raised at that meeting.  Finally, 
the author expresses his deepest thanks to his 
advisor and mentor, Professor Marius Dadarlat, for his support,
encouragement and   the countless hours 
of patiently answering questions.

\newpage

{\large\bf Section 1: \ Preliminaries}

\vspace{10mm}

     In this section we will recall a few facts and present and easy lemma
that will be useful later on.  

{\em Throughout this paper, A and B will
always denote (not necessarily unital) AF algebras}.  Following the 
notation and terminology in [Da], we will let $K_{0}^{+}(A)$ and
$\Gamma(A)$ be the positive cone and scale, respectively, of $K_0(A)$.  
We will say that a group homomorphism, $\theta$: $K_0(A)$ $\rightarrow$ 
$K_0(B)$, is a {\em contraction} if $\theta(\Gamma(A)$) 
$\subset \Gamma(B)$. 
We will say $\theta$ is {\em faithful} if it is contractive and 
$\Gamma(A) \  \cap Ker(\theta$)
 = $\{0\}$.  

 Both of the following two facts  
are essentially contained in Elliot's original classification of AF algebras 
([El]).  Section IV.4 in [Da] has a very nice treatment of this result.
Moreover, the following two facts follow easily from Lemma IV.4.2  and 
the proof of Elliot's Theorem (Theorem IV.4.3) in [Da]. 

\vspace{2mm}

Fact 1.1 \ \ {\em If $\theta$: $K_0(A)$ $\rightarrow$ $K_0(B)$ is a faithful 
group homomorphism  then
there exists a *-monomorphism, $\varphi$: A $\rightarrow$ B, with
$\varphi_{*} = \theta$}.

\vspace{2mm}

Fact 1.2 \ \ {\em If $\varphi , \sigma$: A $\rightarrow$ B are 
*-homomorphisms, 
where A is finite dimensional, then $\varphi_{*} = \sigma_{*} 
\Longleftrightarrow \varphi = Adu \circ \sigma$ for some unitary, u $\in$ 
B.}

\vspace{2mm}

     Recall that when $A$ is unital, $A\times_{\alpha}\Bbb Z$ is defined 
to be the universal
C$^{*}$-algebra generated by $A$ and a unitary $u$, subject to the relation 
$Adu(a) = \alpha(a)$, for all $a \in A$.  We will call $u$ the 
{\em distinguished} unitary in $A\times_{\alpha}\Bbb Z$.  Since one always 
has the freedom to tensor with other AF algebras when proving AF 
embeddability, the following lemma will be useful.

\vspace{2mm}

{\bf Lemma 1.3} \ \ If $A$ and $B$ are unital, $\alpha \in \Aut (A)$ 
and $\beta \in \Aut (B)$ 
are given, then there is a natural embedding

\begin{center} $A\otimes B\times_{\alpha\otimes \beta}\Bbb Z$ 
               ${\Large \hookrightarrow}$
($A\times_{\alpha}\Bbb Z)\otimes (B\times_{\beta}\Bbb Z$)
\end{center}

\vspace{2mm}

{\bf Proof} \ \ Let $u \in A\times_{\alpha}\Bbb Z$, $v \in 
B\times_{\beta}\Bbb Z$, and $w \in 
A\otimes B\times_{\alpha\otimes \beta}\Bbb Z$ be the respective
distinguished unitaries.  Then define a covariant representation by

\begin{center} $w \longmapsto u\otimes v$ \\
               $a\otimes b \longmapsto a\otimes b$
\end{center}

     Let $\varphi: A\otimes B\times_{\alpha\otimes\beta}\Bbb Z 
 \rightarrow (A\times_{\alpha}\Bbb Z)\otimes (B\times_{\beta}
\Bbb Z$) be the induced *-homomorphism and $D$ = range($\varphi$).
  To prove that $\varphi$ is injective (see Thm. 4 in [La]), 
we must provide automorphisms, $\rho(\xi) \in \Aut (D)$, for all $\xi \in
\Bbb C$ with 
$|\xi|$ = 1, such that $\rho(\xi)(a\otimes b) = a\otimes b$ and 
$\rho(\xi)(u\otimes v) = \xi(u\otimes v)$. So, take $\rho_{A}(\xi) \in 
\Aut(A\times_{\alpha}\Bbb Z)$ such that $\rho_{A}(\xi)(a) = a, \ \forall \ 
a \in A,$ and $\rho_{A}(\xi)(u) = \xi u$.  Then let $\rho(\xi) = 
\rho_{A}(\xi) \otimes id_{B\times_{\beta}\Bbb Z}$ and it is easy 
to check that $\rho(\xi)(D) = D$ and $\rho(\xi)|_{D}$ satisfies the 
desired properties. $\Box$

\vspace{7mm}

{\large\bf Section 2: The Rohlin Property}

\vspace{7mm}

     We now state the definition of the Rohlin property in the 
unital case.  The appropriate definition in the non-unital case can be 
found in [EK].

\vspace{4mm}

{\bf Definition 2.1} \ \ If $B$ is  unital  and $\beta \in 
\Aut (B)$ then $\beta$ satisfies the Rohlin property if for every $k \in 
\Bbb N$ there are positive integers $k_{1}, \ldots , k_{m} \geq k$ 
satisfying the following condition: \ For every finite subset, $\cal F$ 
$\subset B$, and every $\epsilon >$ 0, there exist projections 
$e_{i,j}, \  i = 1, \ldots , m, \  j = 0,    \ldots , k_{i} - 1$ in $B$ such 
that

\begin{center}$\sum_{i=1}^{m} \sum_{j=0}^{k_{i} - 1} e_{i,j} = 1_{B}$
\end{center}

\begin{center}$\| \beta(e_{i,j}) - e_{i,j+1} \| < \epsilon$\end{center}
\begin{center}$\| [x,e_{i,j}] \| < \epsilon $\end{center}
for $i = 1, \ldots , m, \ j = 0, \ldots , k_{i} - 1$ and all $x \in \cal F$
(where $e_{i,k_{i}} = e_{i,0}$).  For each integer $i$ above, the 
projections $\{e_{i,j}\}_{j = 0}^{k_i - 1}$ are 
called a {\it Rohlin tower}.

\vspace{4mm}

{\bf Example 2.2} \ \ It follows from Theorem 1.3 in [Ki2] that every UHF 
algebra admits an automorphism with the Rohlin property.  We now give a 
concrete example of such an automorphism which will be used repeatedly. 
Let $M_n$ be the $n\times n$ matrices over $\Bbb C$ 
with canonical matrix units $E_{i,j}^{(n)}$.  Now let $u_n \in M_n$ be the 
unitary matrix such that $Adu_n (E_{i,i}^{(n)}) = 
E_{i + 1, i + 1}^{(n)}$ (with 
addition modulo $n$).  Now consider the Universal UHF algebra 
$\cal U = \otimes_{n \geq 1} M_n$ and the automorphism $\sigma 
= \otimes_{n \geq 1} Adu_n \in \Aut(\cal U)$.  We claim that $\sigma$ 
has the Rohlin property.  It is important to note that the integers 
$k_1, \ldots, k_m$ in the definition of the Rohlin property must be 
fixed and must work for all finite subsets and 
for all choices of $\epsilon$.  We will show that for this particular 
example, one may always take $m = 1$ and $k_1 = k$.  So, let $\cal F 
\subset \cal U$ be a finite subset and $\epsilon > 0$ be given.  Now, 
take  some large $m \in \Bbb N$ such 
that $\cal F$ is nearly contained (within $\epsilon /2$ will suffice) in 
$ \otimes_{n = 1}^{m}M_n \subset \cal U$.  Now take $m^{\prime} > m$ 
such that $k$ divides 
$m^{\prime}$, say $m^{\prime} = sk$.  Then let $e_j = \sum_{t = 0}^{s - 1} 
E_{j + tk, j + tk}^{(m^{\prime})}  \in \cal U$ 
(where we have identified $M_{m^{\prime}}$ with it's 
unital image in $\cal U$), for $1 \leq j \leq k$.  Evidently we have 
$\sigma(e_j) = e_{j + 1}$, 
$\sum_{j = 1}^{k} e_j = 1_{\cal U}$ and the $e_j$'s nearly commute with 
$\cal F$ since they commute with $\otimes_{n = 1}^{m} M_n$ by construction.

     We would also like to point out that $\cal U\times_{\sigma}\Bbb Z$ is 
AF embeddable (and actually embeds back into $\cal U$) 
 by Lemma 2.8 in [Vo].  In fact every crossed product of a UHF algebra 
(by $\Bbb Z$) is AF embeddable since every automorphism of a UHF algebra 
is approximately inner (cf. Theorem 3.6 in [Vo]).  However, this particular 
example is also limit periodic (see Definition 5.1) and thus we do 
not need the full power of Theorem 3.6 in [Vo] to deduce AF embeddability.

\vspace{2mm}

{\bf Remark 2.3} \ \ If $B_0 \subset B$ is a unital, finite dimensional 
subalgebra and 
$\Psi : B \rightarrow B_{0}^{\prime} \cap B$ is a conditional expectation 
then it is readily verified that $\| \Psi(b) - b \|$ is small 
whenever $b$ almost commutes with the set of matrix units from which 
$\Psi$ is constructed.  Thus, if the finite subset, $\cal F$, in the 
definition of the Rohlin property is taken to be the matrix units of 
$B_0$ then we have that the projections, $\{e_{i,j}\}$, are nearly 
contained in $B_{0}^{\prime} \cap B$ and so we can assume without 
loss of generality that $\{e_{i,j}\} \subset B_{0}^{\prime} \cap B$. 
(see Lemma III.3.1 in [Da]).  More generally, if $\cal F$ is contained 
in a finite dimensional subalgebra then the Rohlin towers may be chosen 
to commute with $\cal F$.

\vspace{2mm}

     Throughout this section we will mainly be concerned with unital 
algebras.  However, this is only out of convenience and it should be
noted that all of the results in this section have analogues in the
non-unital case (cf. [EK]).

     We now state a theorem due to Evans and Kishimoto that illustrates 
the usefulness of the Rohlin property when dealing with crossed products. 
Actually, Theorem 4.1 in [EK] is much stronger than the following. 

\vspace{2mm}

{\bf Theorem 2.4}(Evans and Kishimoto) \ \ If $A$ is unital, $\alpha, 
\beta \in \Aut (A)$, 
$\alpha_{*} = \beta_{*}$ on $K_0(A)$, and both $\alpha$ and $\beta$ 
satisfy the Rohlin
property, then $A\times_{\alpha}\Bbb Z \cong A\times_{\beta}\Bbb Z$.

\vspace{2mm}

{\bf Remark 2.5} \ \ The following corollary follows immediately from 
condition (4) of Theorem 0.2 (even without the assumption that $A$ is 
unital). However, it is easily verified that if $\alpha$ has the 
Rohlin property then $\alpha\otimes\beta$ does also.  Thus using 
Theorem 2.4 and Lemma 1.3 we may easily prove homotopy invariance of 
AF embeddability.  

\vspace{2mm}

{\bf Corollary 2.6} \ \ If $A$ is unital and $\alpha, \beta \in \Aut (A)$ 
are homotopic then
$A\times_{\alpha}\Bbb Z$ is AF embeddable $\Longleftrightarrow 
A\times_{\beta}\Bbb Z$ is AF embeddable.

\vspace{2mm}

{\bf Proof} \ \  Let $\cal U$ be the Universal UHF algebra and $\sigma 
\in \Aut(\cal U)$ be the automorphism with the Rohlin property described 
in Example 2.2.  Then $\alpha\otimes \sigma$ is homotopic to 
$\beta\otimes \sigma$.  Hence, from Theorem 2.4 we have 

\begin{center}
  $A\otimes \cal U\times_{\alpha\otimes \sigma}\Bbb Z \cong 
  A\otimes \cal U\times_{\beta\otimes \sigma}\Bbb Z$
\end{center}

     The conclusion now follows from Lemma 1.3 since the crossed product, 
$\cal U\times_{\sigma}\Bbb Z$ is AF embeddable.  $\Box$

\vspace{4mm}

     The following stabilization lemma was essentially proved in [EK] 
(without the assumption of a unit).  Our formulation is somewhat 
different, but 
our proof is just a more detailed version of that given in [EK].

\vspace{2mm}

{\bf Lemma 2.7} \ \ Let $B$ be unital , $\beta \in \Aut (B)$ have the
Rohlin property and $k_{1} \geq k_{2} \geq \ldots \geq k_{m} \geq k$ be 
integers satisfying the definition of the Rohlin property.  If $B_{0} 
\subset B_{1} \subset B_{2}$ are finite dimensional subalgebras of $B$ 
with $1_{B} \in B_{i}, \  i=0,1,2$ and

\begin{center}$\beta^{-j}(B_{i}) \subset B_{i+1}, \ \ $for $ \  i=0,1 \ 
$ and $ 0 \leq j \leq k_{1}$\end{center}
then for every unitary, $u \in B_{2}^{\prime} \cap B$, there exists 
a unitary, $v \in B_{0}^{\prime} \cap B$ such that $\| u - v\beta(v^{*}) 
\| \leq 4/(k - 1)$.
 
\vspace{2mm}

{\bf Proof} \ (cf. [EK], Lemma 3.2) \ Let $u \in B_{2}^{\prime} \cap B$ 
be given.  Then define\\

\begin{center} $\tilde{u}_{0} = 1_B, \ \tilde{u}_{1} = u, \  
\tilde{u}_{j} = 
u\beta(u)\cdots\beta^{j-1}(u),$ \ \ for $j \geq 2$\end{center}
Notice that our hypotheses imply that for $0 \leq j \leq k_{1}$ we have
$\tilde{u}_{j} \in B_{1}^{\prime} \cap B$ and for every unitary $w \in 
B_{1}^{\prime} \cap B$ we have that $\beta^{j}(w) \in B_{0}^{\prime} 
\cap B, \ j = 0, \ldots , k_{1}$.  
     Now, since $B_{1}^{\prime} \cap B$ is a unital AF subalgebra of $B$ 
(cf. Exercise 7.7.5 in [Bl]), 
we can find a unital, finite dimensional subalgebra, $B_{3} \subset 
B_{1}^{\prime} \cap B$, 
such that  there exist unitaries, 
$u_{i} \in B_{3}$ with $\| u_{i} - \tilde{u}_{i} \| 
< \epsilon$, for some small $\epsilon > 0$.  Then let $w^{(t)}_{k_{i}}, i = 1, 
\ldots , m$, 
be paths of unitaries in $ B_{3}$ such that
 
\begin{center} $w^{(0)}_{k_{i}} = 1_{B} \ \ \ \ \ , \ \ \ \ \ 
w^{(1)}_{k_{i}} = u_{k_{i}}$\end{center}

\begin{center} $\| w^{(s)}_{k_{i}} - w^{(t)}_{k_{i}} \| \leq \pi|s - t| \ , 
\ s,t \in [0,1]$\end{center}

     We have already observed that $\beta^{j}(w^{(t)}_{k_{i}}) \in 
B_{0}^{\prime} \cap B$ for $0 \leq j \leq k_1$.  So, (by Remark 2.3) 
take a set of Rohlin 
towers, $\{e_{i,j}\} \subset B_{0}^{\prime} \cap B$.  We may also assume 
without loss of generality that the Rohlin towers approximately commute 
(with error also bounded by $\epsilon$) with the $\tilde{u}_j$'s, 
$\beta(\tilde{u}_j)$'s and 
$\beta^{j}(w^{(1 - \frac{j}{k_i - 1})}_{k_i})$'s below. So we define

\begin{center} $$v^{\prime} = \sum_{i=1}^{m} \sum_{j=0}^{k_{i} - 1} 
\tilde{u}_{j}\beta^{j}
(w^{(1 - \frac{j}{k_{i} - 1})}_{k_{i}})e_{i,j}$$\end{center}

     Now, $v^{\prime}$ is not a unitary.  However, it is easy to verify 
that $\| v^{\prime}v^{\prime*} - 1_{B} \| < \epsilon C$ and 
$\| v^{\prime*}v^{\prime} - 1_B \| < \epsilon C$, where $C$ is a constant 
depending only on $k_1, \ldots, k_m$.  For the 
remainder of the proof, $C$ will always denote a constant depending only 
on $k_1, \ldots, k_m$.  We make no effort to keep track of the best 
constants as this detail is not needed for the proof.  So we 
have that $v^{\prime}$ is invertible and hence the unitary, $v$, in it's polar 
decomposition actually lives in $B_{0}^{\prime} \cap B$ (since 
$v^{\prime} \in B_{0}^{\prime} \cap B$ also).  However, the same estimates 
imply that $v^{\prime}$ is also close to this unitary since $|v^{\prime}|$ 
will be close to $1_B$. More precisely, if we write $v^{\prime} = 
v|v^{\prime}|$ then from elementary spectral theory we have 
$\| v - v^{\prime} \| < 1 - \sqrt{1 + \epsilon C}$.

     Notice that the definition of the Rohlin property implies the 
following estimates.

\begin{center}
  $\| \beta(e_{i,j})e_{i,j + 1} - e_{i, j + 1} \| \leq \epsilon$
\end{center}
\begin{center}
  $\| \beta(e_{i,j})e_{k,l} \| \leq \epsilon$, \ unless $k = i$ and $l = j + 1$
\end{center}
Thus a straightforward calculation yields
\begin{center}
  $$\| v^{\prime}\beta(v^{\prime*}) - 
  \sum_{i = 1}^{m}\sum_{j = 0}^{k_i - 1}\tilde{u}_{j}\beta^{j}(w_{k_i}^{(1 - 
  \frac{j}{k_i - 1})})e_{i,j}\beta^{j}(w_{k_i}^{(1 - \frac{j - 1}{k_i - 1})*})
\beta(\tilde{u}_{j - 1}^{*}) \| \leq \epsilon C$$
\end{center}
where  for each 
$i$, the $j = 0$ term in the sum above is actually given by
\begin{center}
 $\tilde{u}_0 u_{k_i}e_{i,0}\beta^{k_i}(w_{k_i}^{(1 - 
 \frac{k_i - 1}{k_i - 1})*})\beta(\tilde{u}_{k_i - 1}^{*}) = 
 u_{k_i}e_{i,0}\beta(\tilde{u}_{k_i - 1}^{*}) \approx \tilde{u}_{k_i} 
 e_{i,0}\beta(\tilde{u}^{*}_{k_i - 1})$
\end{center} 
\vspace{1mm}

    The following observations, will be needed to get the desired estimate.

\begin{enumerate}
\item[a)] $\| e_{i,j}\beta^{j}(w_{k_i}^{(1 - \frac{j - 1}{k_i - 1})*}) 
           \beta(\tilde{u}_{j - 1}^{*}) - 
           \beta^{j}(w_{k_i}^{(1 - \frac{j - 1}{k_i - 1})*})
           \beta(\tilde{u}_{j - 1}^{*})e_{i,j} \| \leq 2\epsilon$
\item[b)] $\| \beta^{j}(w_{k_i}^{(1 - \frac{j}{k_i - 1})})
              \beta^{j}(w_{k_i}^{(1 - \frac{j -1}{k_i - 1})*}) - 1_{B} \| 
              \leq \frac{\pi}{k_i - 1}$
\end{enumerate}
\begin{enumerate}
\item[c)] $\tilde{u}_j \beta(\tilde{u}_{j - 1}^{*}) = u$ \ for $j \geq 
          1 $
\end{enumerate}        

     Using the estimates above, we have the following approximations.

\begin{eqnarray*}
  v^{\prime}\beta(v^{\prime*}) & \stackrel{\epsilon C}{\approx} &
 \sum_{i = 1}^{m}\sum_{j = 0}^{k_i - 1}\tilde{u}_{j}\beta^{j}(w_{k_i}^{(1 - 
  \frac{j}{k_i - 1})})e_{i,j}\beta^{j}(w_{k_i}^{(1 - \frac{j - 1}{k_i - 1})*})
\beta(\tilde{u}_{j - 1}^{*}) \\ & \stackrel{\epsilon C}{\approx} &
 \sum_{i = 1}^{m}\sum_{j = 0}^{k_i - 1}\tilde{u}_{j}\beta^{j}(w_{k_i}^{(1 - 
  \frac{j}{k_i - 1})})\beta^{j}(w_{k_i}^{(1 - \frac{j - 1}{k_i - 1})*})
\beta(\tilde{u}_{j - 1}^{*})e_{i,j} \ \ \ \ by \ \ (a) \\ & 
\stackrel{\frac{\pi}{k - 1} + \epsilon C}{\approx} &
\sum_{i = 1}^{m}\sum_{j = 0}^{k_i - 1}\tilde{u}_{j}
\beta(\tilde{u}_{j - 1}^{*})e_{i,j}  \hspace{50mm} 
by \ \ (b) \\ & \stackrel{\epsilon C}{\approx} &
u \hspace{80mm}  by \ \ (c).        
\end{eqnarray*}                         
                                         
     Recall that the only thing keeping the last approximation from being 
an equality is the fact that for each $i$, the $j = 0$ term above is actually 
given by $u_{k_i} \beta (\tilde{u}_{k_i - 1}^{*})e_{i,0}  
\stackrel{\epsilon}{\approx} \tilde{u}_{k_i} 
\beta(\tilde{u}_{k_i - 1}^{*})e_{i,0} = ue_{i,0}$. Finally, 
combining all of these 
estimates  we have that  $\| u - v^{\prime}\beta(v^{\prime*}) \| 
< \frac{\pi}{k - 1} + \epsilon C$ and hence $\| u - v\beta(v^{*}) \| < 
\frac{\pi}{k - 1} + \epsilon C + 2(1 - \sqrt{1 + \epsilon C})$,
  where we have the liberty to take 
$\epsilon$ as small as we like.  $\Box$

\vspace{4mm}

     We now turn to the main tool of this paper.  The following proposition
is the key to Proposition 3.1, which in turn is the key to the rest of the 
embedding results presented here.

\vspace{4mm}

{\bf Proposition 2.8} \ \ Let $A, B$ be unital, $\alpha \in \Aut (A)$,
$\beta \in \Aut (B)$ and $\varphi: A \rightarrow B$ a unital, *-monomorphism.
Assume also that $\beta$ satisfies the Rohlin property and 
$\beta_{*} \circ \varphi_{*} = \varphi_{*} \circ \alpha_{*}$, i.e. 
that we have commutativity in the diagram
\vspace{2mm}

\begin{center} $K_{0}(A) \stackrel{\varphi_{*}}{\longrightarrow} K_{0}(B)$ 
\end{center}

\begin{center}  $\alpha_{*}$ {\large $\downarrow$} \ \ \ \ \ \ \ \ \ \ \ \
{\large $\downarrow$} $\beta_{*}$ \end{center}

\begin{center} $K_{0}(A) \stackrel{\varphi_{*}}{\longrightarrow} K_{0}(B)$
\end{center}
Then there exist unitaries, $v \in A, u \in B$, and a unital 
*-monomorphism,\\
$\varphi^{\prime}: A \rightarrow B$, with $\varphi_{*}^{\prime} = 
\varphi_{*}$ and commutativity in 

\begin{center} $A \stackrel{\varphi^{\prime}}{\longrightarrow} B$

               {\small $Adv \circ \alpha$} {\Large $\downarrow$} \ \ \ \ \ \  
{\Large $\downarrow$} {\small $Adu \circ \beta$}

               $A \stackrel{\varphi^{\prime}}{\longrightarrow} B$
\end{center}

\vspace{2mm}

{\bf Proof} \ \ Let $\{ A_{i} \}_{i=0}^{\infty}$ be an increasing nest of
finite dimensional subalgebras whose union is dense in $A$ and $A_{0} 
= \Bbb C 1_{A}$.
Then $\{ \alpha(A_{i}) \}$ is a nest with the same properties and thus we can
find a unitary $v \in A$ (with $v$ as close to $1_{A}$ as we like) 
such that (see Theorem III.3.5 in [Da])

\begin{center} $Adv\circ\alpha(\bigcup A_{i}) = \bigcup A_{i}$ \end{center}
To ease our notation somewhat, we will henceforth denote $Adv\circ\alpha$ also
by $\alpha$ (and just remember that the $v$ in the statement of the proposition
has already been chosen).

     Now, let $m_i = 2^{i + 3} + 1$ and applying the definition of the 
Rohlin property to $m_i$ we get a finite set of integers, 
$\{k_{1}^{(i)}, \ldots, k_{s_i}^{(i)}\}$, with $k_{t}^{(i)} \geq m_i$ for 
$1 \leq t \leq s_i$.  Then let $m_{i}^{\prime} = max\{k_{1}^{(i)}, 
\ldots, k_{s_i}^{(i)}\}$.

    Since we have arranged that $\bigcup\alpha(A_{i}) = \bigcup
A_{i}$, by passing to a subsequence we may further assume that for $i \geq 2$
we have that

\begin{center} $\alpha^{-j}(A_{i - 1}) \subset A_{i} , \ 0
\leq j \leq m_{i}^{\prime}$ \end{center}
\begin{center} $\alpha^{-j}(A_{i - 2}) \subset A_{i - 1} , \ 0 
\leq j \leq m_{i}^{\prime}$ \end{center}

     We will now inductively construct sequences of unitaries, $u_{i} \in B$,
automorphisms, $\beta_{i} \in \Aut (B)$, and unital *-monomorphisms, 
$\varphi_{i}: A \rightarrow B$, with the following properties

\begin{center}
 \begin{enumerate}
 \item $\| u_{i} - 1_{B} \| \leq 1/2^{i}$
 \item $\beta_{i + 1} = {\rm Adu}_{i + 1}\circ\beta_{i}$
 \item $\varphi_{i*} = \varphi_{*}$
 \item $\varphi_{i + 1}|_{\alpha(A_{i - 2})} = 
        \varphi_{i}|_{\alpha(A_{i - 2})} \ $ where $ \  \Bbb C 1_{A} = A_{0} = 
         A_{-1} = A_{-2}$
 \item $\beta_{i}\circ\varphi_{i}|_{A_{i}} = \varphi_{i}\circ\alpha|_{A_{i}}$  
 \end{enumerate}
\end{center}
where $u_{0} = 1_{B}, \  \varphi_{0} = \varphi$ and $\beta_{0} = \beta$. 
Note that each $\beta_i$ will also have the Rohlin property.  So, 
assume that we have found the desired $u_{i}, \  \beta_{i}, \ 
 \varphi_{i}$ for 
$0 \leq i \leq n$ and we will show how  to construct $u_{n + 1}, \ 
 \beta_{n + 1}$ and
  $ \varphi_{n + 1}$.  By construction, we have $\beta_{n*}\circ\varphi_{n*} = 
\varphi_{n*}\circ\alpha_{*}$ and thus restricting to $A_{n + 1}$ we have that
$\beta_{n*}\circ\varphi_{n*}$ and $\varphi_{n*}\circ\alpha_{*}$ agree as maps
from $K_{0}(A_{n + 1}) \rightarrow K_{0}(B)$.  Thus, (by Fact 1.2) 
we can find a unitary 
$w_{n + 1} \in B$ such that

\begin{center} $Adw_{n + 1}\circ\beta_{n}\circ\varphi_{n}|_{A_{n + 1}} = 
\varphi_{n}\circ\alpha|_{A_{n + 1}}$ \end{center}
But, by assumption, we already have

\begin{center} $\beta_{n}\circ\varphi_{n}|_{A_{n}} = 
                \varphi_{n}\circ\alpha|_{A_{n}}$\end{center}
and hence $w_{n + 1} \in \varphi_{n}(\alpha(A_{n}))^{\prime} \cap B$.  Now, 
since we have control of the iterates of $A_{n - 2}$ and $A_{n - 1}$ (under 
$\alpha$) we claim that

\begin{center} $\beta_{n}^{-j}(\varphi_{n}(\alpha(A_{n - 1}))) \subset 
\varphi_{n}(\alpha(A_{n})) \  \ , \ \ 0 \leq j \leq m_{n}^{\prime}$
 \end{center}
\begin{center} $\beta_{n}^{-j}(\varphi_{n}(\alpha(A_{n - 2}))) \subset
\varphi_{n}(\alpha(A_{n - 1})) \ \ , \ \ 0 \leq j \leq m_{n}^{\prime}$
 \end{center}
To see this we first note that since $\alpha^{-1}(A_{n - 1}) 
\subset A_n$,  $(5)$ implies  
\begin{center}
  $\beta^{-1}_{n}(\varphi_{n}(\alpha(A_{n -1}))) = \varphi_{n}(A_{n - 1})
  \subset \varphi(\alpha(A_n))$
\end{center}
Similarly we have  
\begin{eqnarray*}
  \beta^{-2}_n(\varphi_n(\alpha(A_{n - 1}))) & = &
  \beta^{-1}(\beta^{-1}(\varphi_n(\alpha(A_{n - 1}))))\\
  & = & 
  \beta^{-1}(\varphi_n(A_{n - 1}))\\
  & = &
  \beta^{-1}(\varphi_n(\alpha(\alpha^{-1}(A_{n - 1}))))\\
  & = &
  \varphi_n(\alpha^{-1}(A_{n - 1})) \ \ \ \ ( since \ \ 
  \alpha^{-1}(A_{n - 1}) \subset A_n)
\end{eqnarray*}
but, $\alpha^{-2}(A_{n - 1}) \subset A_n$ implies that 
$\varphi_n(\alpha^{-1}(A_{n - 1})) \subset \varphi_n(\alpha(A_n))$.  By 
repeating the above argument, one can show that 
$\beta_{n}^{-j}(\varphi_n(\alpha(A_{n - 1}))) = 
\varphi_n(\alpha^{-j + 1}(A_{n - 1}))$ for $0 \leq j \leq m_{n}^{\prime}$. 
But then $\alpha^{-j}(A_{n -1}) \subset A_n$ implies that 
$\varphi_n(\alpha^{-j + 1}(A_{n -1})) \subset \varphi_n(\alpha(A_n))$, 
for $0 \leq j \leq m_{n}^{\prime}$.  The same argument works for the 
iterates of $\varphi_n(\alpha(A_{n - 2}))$ as well.

Thus, by Lemma 2.7, we may take a unitary, $v_{n + 1} \in 
\varphi_{n}(\alpha(A_{n - 2}))^{\prime} \cap B$ such that $\| w_{n + 1} - 
v_{n + 1}\beta_{n}(v_{n + 1}^{*}) \| \leq \frac{4}{m_{n} - 1} = 
\frac{1}{2^{n + 1}}$.
  So, we define
\begin{center}
$\varphi_{n + 1} = Adv_{n + 1}^{*}\circ\varphi_{n}$,

 $ u_{n + 1} = v_{n + 1}^{*}w_{n + 1}\beta_{n}(v_{n + 1})$
 
$ \beta_{n + 1} = Adu_{n + 1}\circ\beta_{n}$
\end{center}
It is now easy to check that we have 
satisfied all of the required properties.

     The proof of the proposition is now complete as it is easy to see that
$u = \lim_{n} u_{n}\cdots u_{1}$ is a well defined unitary in $B$.  Also, it
is easy to see (by condition 4) that the $\varphi_{i}$'s converge (in the
point-norm topology) to a unital *-monomorphism, $\varphi^{\prime}$, and
clearly $Adu\circ\beta\circ\varphi^{\prime} = \varphi^{\prime}\circ\alpha$.  
$\Box$   

\vspace{4mm}
 
{\bf Remark 2.9} \ \ It would be of independent interest to know if Proposition
2.6 can be proved without the assumption that $\beta$ satisfy the Rohlin
property.  However, the author was unable to either prove or provide a 
counterexample to such a claim.

\vspace{4mm}

{\bf Corollary 2.10} \ \ Under the assumptions of Proposition 2.8 we also have
that $A\times_{\alpha}\Bbb Z$ embeds into $B\times_{\beta}\Bbb Z$.

\vspace{2mm}

{\bf Proof} \ \  Essentially what we have done in Proposition 2.8 is embedded 
$A$ into $B$ in such a way that the automorphism $Adv \circ \alpha$ 
extends to an automorphism (namely, $Adu \circ \beta$) of $B$.  
In general it is
true that if $C, \ D$ are C$^{*}$-algebras with $C \subset D$, 
$\gamma \in \Aut (D)$ and $\gamma (C) = C$ then there is a natural inclusion

\begin{center}
  $C\times_{\gamma|_{C}}\Bbb Z \hookrightarrow D\times_{\gamma}\Bbb Z.$
\end{center}

     Thus the conclusion follows from the isomorphisms
\begin{center}
   $A\times_{\alpha}\Bbb Z \cong A\times_{Adv\circ\alpha}\Bbb Z, \ 
 \ B\times_{\beta}\Bbb Z \cong B\times_{Adu\circ\beta}\Bbb Z$ 
\end{center}
and  the fact that $Adu \circ \beta$ is an extension of $Adv \circ \alpha$ 
under the embedding $\varphi^{\prime}$. $\Box$  

\vspace{7mm}

{\large\bf Section 3: Characterizing AF Embeddability}

\vspace{7mm}

     If $\varphi: A \rightarrow B$ is a *-homomorphism, we will denote 
by $\tilde{A}$ and $\tilde{B}$ the C$^{*}$-algebras obtained by adjoining 
a (possibly new) unit and we will let $\tilde{\varphi}$ be the (unique) 
unital extension of $\varphi$.

\vspace{4mm}

{\bf Proposition 3.1} \ \ Let $A$ and $B$ be given (not necessarily unital), 
$\alpha \in \Aut (A), \beta \in \Aut (B)$ and $\varphi: A \rightarrow B$ be a 
*-monomorphism.   Further assume that $\beta_{*}\circ\varphi_{*} = 
\varphi_{*}\circ\alpha_{*}$ and that 
$\tilde{B}\times_{\tilde{\beta}}\Bbb Z$ is AF embeddable.
Then $A\times_{\alpha}\Bbb Z$ is also AF embeddable.

\vspace{2mm}

{\bf Proof} \ \ Note that after unitizing everything we still have 
$\tilde{\beta_{*}}\circ\tilde{\varphi_{*}} = 
\tilde{\varphi_{*}}\circ\tilde{\alpha_{*}}$.  
  Let $\cal U$ be the universal 
UHF algebra and $\sigma \in \Aut (\cal U$) be the automorphism 
(with the Rohlin property) defined in Example 2.2.
  Now we let 
$\tilde{\beta}^{\prime} 
= \tilde{\beta}\otimes \sigma$ and let 
$\tilde{\varphi}^{\prime}: A \rightarrow 
\tilde{B}\otimes \cal U$ be given by $\tilde{\varphi}^{\prime}(a) 
= \tilde{\varphi}(a)\otimes 
 1_{\cal U}$.  Then $\tilde{\beta}^{\prime}$ satisfies the Rohlin 
property and it is
easy to check that we still have commutativity at the level of K-theory, i.e. 
$\tilde{\beta}_{*}^{\prime}\circ\tilde{\varphi}_{*}^{\prime} = 
\tilde{\varphi}_{*}^{\prime}\circ\tilde{\alpha}_{*}$. Now, it is always 
true that 

\begin{center}
  $A\times_{\alpha}\Bbb Z$ $\hookrightarrow$ 
    $\tilde{\rm A}\times_{\tilde{\alpha}}\Bbb Z$ 
\end{center}
and by Corollary 2.10 we have that

\begin{center} $\tilde{A}\times_{\tilde{\alpha}}\Bbb Z$ 
   $\hookrightarrow$ 
   $\tilde{B}\otimes \cal U\times_{\tilde{\beta}^{\prime}}\Bbb Z$.
\end{center}
But then by Lemma 1.3 we also have that

\begin{center} 
  $\tilde{B}\otimes \cal U\times_{\tilde{\beta}^{\prime}}\Bbb Z$ 
  $\hookrightarrow$ 
  ($\tilde{B}\times_{\tilde{\beta}}\Bbb Z)\otimes
   (\cal U\times_{\sigma}\Bbb Z$) 
\end{center}
where each of the algebras on the right hand side are AF embeddable. $\Box$

\vspace{2mm}

{\bf Corollary 3.2} \ \ If $\alpha \in \Aut (A)$ and $\theta : K_0 (A) 
\rightarrow G$ is a positive group homomorphism where 
$G$ is a dimension group (cf. [Ef]),   
$Ker (\theta) \cap \Gamma(A) = \{0\}$ and $H_{\alpha} \subset Ker(\theta)$ 
 then $A\times_{\alpha}\Bbb Z$ is AF embeddable.

\vspace{2mm}

{\bf Proof} \ \ Since there exists 
an AF algebra, $B$, with $K_0 (B) = G$ and $\Gamma(B) = G^{+}$,  we have 
that $\theta$ is a faithful homomorphism and thus by Fact 1.1 
 we may find a *-monomorphism, 
$\varphi : \ A \rightarrow B$, with $\varphi_{*} = \theta$.  
But we may now take $\beta 
\in \Aut (B)$ to be the identity and the conclusion  follows from 
Proposition 3.1. $\Box$

\vspace{2mm}
 
{\bf Remark 3.3} \ \ The hypotheses of Corollary 3.2 are easily seen to 
be necessary as well.  The proof of this necessity is essentially contained 
in the introduction. 

\vspace{2mm}

     We nearly have the necessary tools to prove Theorem 0.2 now.  However, 
for the implication (4) $\Rightarrow$ (1) we need to provide a candidate 
AF algebra for the desired embedding of $A\times_{\alpha}\Bbb Z$.  The 
following key lemma due to Spielberg will give us our candidate.  We will 
not prove this lemma (see Lemma 1.14 in [Sp]) but would like to point 
out that it depends on the Effros-Handelman-Shen Theorem ([EHS]). 

\vspace{2mm}

{\bf Lemma 3.4}(Spielberg) \ If $G$ is a dimension group and $H$ is a subgroup 
of $G$ with $H \cap G^{+} = \{0\}$, then there exists a positive group 
homomorphism $\theta: G \rightarrow G^{\prime}$ (where $G^{\prime}$ is also 
a dimension group) with 
\begin{enumerate}
 \item $H \subset Ker(\theta)$
 \item $Ker(\theta) \cap G^{+} = \{0\}$
\end{enumerate}

\vspace{2mm}

{\bf Proof of Theorem 0.2} 

$(1) \Rightarrow (2)$ is obvious (whether or not $A$ is unital).

$(2) \Rightarrow (3)$ If $A$ is unital then so is $A\times_{\alpha}\Bbb Z$ 
and thus quasidiagonality implies that every isometry in 
$A\times_{\alpha}\Bbb Z$ is actually a unitary.  But this implies 
that the identity of $A\times_{\alpha}\Bbb Z$ is not equivalent to any 
proper subprojection and hence $A\times_{\alpha}\Bbb Z$ is finite.  Since 
matrix algebras over quasidiagonal algebras are again quasidiagonal, the 
same argument shows that $A\times_{\alpha}\Bbb Z$ is stably finite.

     If $A$ is non-unital and $A\times_{\alpha}\Bbb Z$ is quasidiagonal then 
so is $\widetilde{A\times_{\alpha}\Bbb Z}$.  Thus by the above 
argument we have that $\widetilde{A\times_{\alpha}\Bbb Z}$ is stably finite 
and thus $A\times_{\alpha}\Bbb Z$ inherits this property also.

$(3) \Rightarrow (4)$ We have already shown this in the introduction when 
$A$ is unital. 
However, the following is a very elementary argument which does not 
depend on the Pimsner-Voiculescu six term exact sequence and which 
actually holds for much general algebras than just AF algebras.

     Since $A$ has a countable approximate unit consisting 
of projections, we do not need to pass to the unitization of $A$ when 
computing $K_0 (A)$ (see Proposition 5.5.5 in [Bl]).  Now, assume 
that $H_{\alpha} \cap K_{0}^{+}(A) 
\supsetneqq \{0\}$, i.e. there is some $x \in K_0 (A)$ and some projection, 
$r \in M_n (A)$ such that $\alpha_{*} (x) - x = [r] \neq 0$.  Now, write 
$x = [p] - [q]$, where (without loss of generality) $p, \ q \in M_n (A)$ 
are projections.  We will not keep track of the sizes of matrices that we 
are dealing with since this is not important for our argument.  So, 
rewriting the equation $\alpha_{*} (x) - x = [r]$ we get 

\begin{center}
  $[\alpha (p)] + [q] = [p] + [\alpha (q)] + [r]$
\end{center} 

     Now, since AF algebras enjoy cancellation, we can find a partial 
isometry (in the matrices over A), $v$, with support projection 
$diag(p, \ \alpha(q), \ r)$ and range projection $diag(\alpha (p), \ q, \ 0)$. 
If we now move to $A\times_{\alpha}\Bbb Z$ and let $u$ be the distinguished 
unitary (i.e. $uau^{*} = \alpha (a), \forall \ a \in A$) then it is easy to 
verify that $diag(pu^{*}, \ uq, \ 0)$ is a partial isometry with support 
projection $diag(\alpha (p), \ q, \ 0)$ and range projection 
$diag(p, \ \alpha (q), \ 0)$.  Thus multiplying these two partial isometries 
we get a partial isometry (in the matrices over $A\times_{\alpha}\Bbb Z$) 
from $diag(p, \ \alpha (q), \ r)$ to the proper subprojection 
$diag(p, \ \alpha(q), \ 0)$ and thus $A\times_{\alpha}\Bbb Z$ is not 
stably finite.

$(4) \Rightarrow (1)$ Letting $H = H_{\alpha}$ in  Lemma 3.4, this 
implication now follows from Corollary 3.2. $\Box$

\vspace{2mm}

{\bf Remark 3.5} \ \ Note that the implications $(1) \Rightarrow (2) 
\Rightarrow (3)$ hold for any C$^{*}$-algebra, $A$, while the only 
property of AF algebras that we really used in the implication 
$(3) \Rightarrow (4)$ was the fact that one need not pass to the 
unitization  of $A$ when computing $K_0(A)$ (i.e. cancellation is 
not necessary in the argument above).

\vspace{7mm}

{\large\bf Section 4: Applications}

\vspace{7mm}

     We now introduce a simple K-theoretical condition which is sufficient
to ensure AF embeddability of the corresponding crossed product.  This 
condition gives a generalization of all the embedding theorems in [Vo] as
all of the automorphisms considered there satisfy the following.

\vspace{2mm}

{\bf Definition 4.1} \ \ If $\alpha \in \Aut (A)$ then 
$\alpha$ satisfies the {\it finite orbit property}, (FOP), if for every $x 
\in K_{0}(A)$ there exists an integer, $0 \neq n \in \Bbb N$, such
that $\alpha_{*}^{n}(x) = x.$
          
\vspace{2mm}

{\bf Theorem 4.2} \ \ If A (not necessarily unital) is given, $\alpha 
\in \Aut (A)$ and $\alpha$ satisfies  (FOP) then $A\times_{\alpha}\Bbb Z$ is
AF embeddable.

{\bf Proof} \ \ We will show the contrapositive.  If 
$A\times_{\alpha}\Bbb Z$ is not AF embeddable then $H_{\alpha} \cap 
K_{0}^{+} (A) \neq \{0\}$.  So, let 
$x \in K_0(A)$ be chosen so that $0 \neq \alpha_{*}(x) - x \in K_{0}^{+}(A)$. 

     Note that since $\alpha_{*}$ is an isomorphism and $K_{0}^{+}(A) \cap 
(-K_{0}^{+}(A)) = 0$ we have the following two facts.

\vspace{1mm}

     $i) If \ 0 \neq [p] \in K_{0}^{+}(A) \  then \ 0 \neq \alpha_{*}([p]) 
= [\alpha(p)] \in K_{0}^{+}(A)$ 

\vspace{1mm}

     $ii) If \ 0 \neq [p], [q] \in K_{0}^{+}(A)  \  then \ 0 \neq [p] + [q] 
 \in K_{0}^{+}(A)$

\vspace{1mm}

     Thus we have that $0 \neq \alpha_{*}^{2}(x) - \alpha_{*}(x)$ and hence 
$0 \neq (\alpha_{*}^{2}(x) - \alpha_{*}(x)) + (\alpha_{*}(x) - x) =
\alpha_{*}^{2}(x) - x \in K_{0}^{+}(A)$.  Arguing similarly we have that 
$0 \neq \alpha_{*}^{j}(x) - x$, for all nonzero $j \in \Bbb N$.  Thus $\alpha$ 
does not satisfy (FOP). $\Box$

\vspace{2mm}

{\bf Remark 4.3} \ \ Recall that Theorem 3.6 in [Vo] states that if 
there exists a nonzero integer, $n$, such that $\alpha_{*}^{n} = id_{A*}$ 
then the corresponding crossed product is AF embeddable.  The following 
example shows that Theorem 4.2 is a generalization of this result.
  Let $X = \{1, 1/2, 1/3, 1/4, \ldots\} \cup \{0\}$ and 
let $\varphi: X \rightarrow X$ be the homeomorphism which leaves 0 fixed, 
interchanges 1 and 1/2, cyclicly permutes 1/3, 1/4 and 1/5, cyclicly 
permutes 1/6, 1/7, 1/8 and 1/9 and so on (with increasing lengths of 
cycles).  Then it is easy to see that $\varphi$ satisfies (FOP), but 
there is no power of $\varphi$ which is approximately inner since the 
only inner automorphism is the identity map.

\vspace{2mm}

     It is not hard to show that (FOP) is not a necessary condition for 
AF embeddability (see Remark 4.8).  However, automorphisms satisfying 
(FOP) do admit a nice characterization.  It should be 
pointed out that the following proposition will not be needed in the 
remainder of this paper.

     Note that $\alpha$ satisfies (FOP) if and only if $\tilde{\alpha}$ 
satisfies (FOP) and thus there is no harm in passing to unitizations 
when dealing with this property.

     Abusing notation a little, we will  also denote by 
$\alpha_{*}$ the image of an automorphism in the quotient group 
$\Aut (A)_{*} = \Aut (A)/\overline{Inn(A)}$.  We now  
show that  (FOP) characterizes all the well behaved automorphisms 
in $\Aut (A)_{*}$.

\vspace{4mm}

{\bf Proposition 4.4} \ \ Let $A$ be  unital and $\alpha \in 
\Aut (A)$.
Then the following are equivalent:

\begin{enumerate}
 \item $\alpha$ satisfies (FOP)
 \item There exists a sequence, $n_{k} > 0$, such that $d(\alpha^{n_{k}}, 
       Inn(A)) \rightarrow$ 0, as $k \rightarrow \infty$ (where $d$ is any
       metric on $\Aut (A)$ which induces the point-norm topology).
 \item The sequence $\{\alpha_{*}^{n}\}_{n \geq 0}$ has a convergent 
       subsequence in $\Aut (A)_{*}$.
\end{enumerate}

{\bf Proof} \ \ (1) $\Rightarrow$ (2) \ \ Let $\{A_{i}\}$ be any increasing
nest of finite dimensional subalgebras whose union is dense in $A$.  It 
suffices to show that for each $i \in \Bbb N$ we can find  integers $n_{i}$
(with $n_{i} <  n_{i + 1}$) and unitaries $u_{i} \in A$ such that 
$\alpha^{n_{i}}|_{A_{i}} = Adu_{i}|_{A_{i}}$.  So let $e_{1}, \ldots , 
e_{m_{i}}$ be minimal projections in $A_{i}$ such that $\{ [e_{1}], \ldots, 
[e_{m_{i}}] \}$ generate $K_{0}(A_{i})$.  Now, take an integer $n_{i}$ such 
that $\alpha^{n_{i}}_{*}$([$e_{j}$]) = [$e_{j}$], for $1 \leq j \leq m_{i}$
(note that $n_{i}$ may be chosen larger than any specified number).  Thus, the
restriction of $\alpha^{n_{i}}$ to $A_{i}$ agrees on K-theory with the natural
inclusion map $A_{i} \hookrightarrow A$.  Thus we may find the desired unitary
, $u_{i} \in A$. 

     (2) $\Rightarrow$ (3) \ \ Clearly $\alpha_{*}^{n_{i}} \rightarrow id_{A*}$

     (3) $\Rightarrow$ (4) \ \ Assume there exists a subsequence $\{m_{k}\}$
and an element $\beta_{*} \in \Aut (A)_{*}$ such that $\alpha_{*}^{m_{k}} 
\rightarrow \beta_{*}$.  Then defining $n_{k} = m_{k} - m_{k - 1}$ we have that

\begin{center} $\alpha_{*}^{n_{k}} = \alpha_{*}^{m_{k}} \circ 
\alpha_{*}^{-m_{k - 1}} \rightarrow \beta_{*} \circ \beta_{*}^{-1} = id_{A*}$
\end{center}

     Since $\Gamma(A)$ generates $K_0(A)$, it 
suffices to check  (FOP) on this set.  So, let $p \in A$ be any projection.
Then let $\{a_{i}\}$ be any sequence which is dense in the unit ball of A, with
$a_{1}$ = p.  Now recall that the metric on Aut($A$) is defined as follows

\begin{center}$$ d(\sigma,\gamma) = \sum_{n = 1}^{\infty} \frac{1}{2^{n}} 
\| \sigma(a_{i}) - \gamma(a_{i}) \| , \ \ {\rm for}  \  \sigma, \gamma \in 
 \Aut (A)$$ \end{center}
where choosing a different sequence just gives an equivalent metric (cf. 
[Ar]).  Now, let $U_{1/2}$ = $\{ \sigma \in \Aut (A): \text{there exists} 
 \ \gamma \in Inn(A) \text{ such that} \ d(\sigma,\gamma) < 1/2\}$.  
Note that if $\sigma \in U_{1/2}$ then $\sigma_{*}([p]) = [p]$, since 
$a_{1} = p$ and $\| \sigma(p) - \gamma(p) \| < 1$ for some $\gamma \in 
Inn(A)$.

     Now, recall that by definition of the topology on $\Aut (A)_{*}$ (cf. 
[Po]) we have that $U_{1/2 *} = \{ \sigma_{*} : \sigma \in
U_{1/2} \}$ is an open set in $\Aut (A)_{*}$ (obviously containing $id_{A*}$).
Thus, there exists a $k \in \Bbb N$ large enough that $\alpha_{*}^{n_{k}} \in
U_{1/2 *}$.  Thus, there exists $\sigma \in U_{1/2}$ with 
$\alpha_{*}^{n_{k}} = \sigma_{*}$ and hence $\alpha_{*}^{n_{k}}([p]) = 
\sigma_{*}([p]) = [p]$. $\Box$

\vspace{4mm}

 We will now show how to recover (the AF case of) a result of 
Pimsner on the AF 
embeddability of crossed products of commutative C$^{*}$-algebras. 
(See [Pi])

\vspace{2mm}

{\bf Definition 4.5} If $X$ is a compact, metrizable space and $\varphi$ is 
a homeomorphism of $X$.  Then $x \in X$ is called {\it pseudo-nonwandering} if 
for every $\epsilon > 0$, there exists a set of points $x_0, \ldots, 
x_{n + 1}$ 
with $x = x_0 = x_{n + 1}$ and $d(\varphi(x_i),x_{i + 1}) < \epsilon$ 
for $0 \leq i \leq n$.  
The set of all pseudo-nonwandering points will 
be denoted by $X(\varphi)$.

\vspace{2mm}

{\bf Lemma 4.6}(Pimsner)  If $\varphi: X \rightarrow X$ is a homeomorphism 
of the compact metrizable space, $X$, then $x \notin X(\varphi) 
\Leftrightarrow$ there exists an open set, $U \subset X$ such that 
$\varphi(\overline{U}) \subset U$ and 
$x \in U\backslash\varphi(\overline{U})$.

\vspace{2mm}

     As we are interested in the case when $C(X)$ is AF, we will assume 
from now on that $X$ is also totally disconnected.  Then we can find a 
basis of clopen sets, say $\{V_n\}$.  We now claim that in this case, the 
open set $U$ in Lemma 4.2 can be taken to be clopen.  To see this, we 
assume that there exists an open set $U \subset X$ satisfying the two 
properties of the lemma.  Now, since $\varphi(\overline{U})$ is a 
compact subset of $U$, we can find a finite subset of the clopen 
basis, say $V_{n_1}, \ldots, V_{n_k}$, such that 

\begin{center}
   $x \notin \cup_{i = 1}^{k}V_{n_i}$\\
\vspace{1mm}
   $\varphi(\overline{U}) \subset \cup_{i = 1}^{k}V_{n_i} \subset U$
\end{center}

     Then letting $V = \varphi^{-1}(\cup_{i = 1}^{k}V_{n_i})$, we get that 
$V$ is a clopen set with the same properties as $U$. 

    Recall that when $X$ is totally disconnected, $K_0 (C(X)) = 
C(X, \Bbb Z)$, the group of continuous functions from $X$ to $\Bbb Z$, 
with positive cone given by the nonnegative functions.  

\vspace{2mm} 

{\bf Theorem 4.7}(Pimsner) \ If $\varphi: X \rightarrow X$ is a 
homeomorphism of the compact, totally disconnected metric space $X$, 
then $C(X)\times_{\varphi}\Bbb Z$ is AF embeddable if and only if 
$X(\varphi) = X$, where the $\varphi$ appearing in the crossed 
product denotes the corresponding automorphism of $C(X)$, i.e. 
$\varphi(f) = f \circ \varphi^{-1}$.

\vspace{2mm}

{\bf Proof} $(\Rightarrow)$  We prove the contrapositive.  So assume 
that $X(\varphi) \neq X$.  Then by Lemma 4.6 and the discussion which follows 
, we can find a clopen set, $V$, such that $\varphi(V)$ is properly 
contained in $V$.  So, if $P \in C(X)$ is the projection with support 
$V$, then we have that $P \circ \varphi^{-1}$ is a projection in $C(X)$ 
which is dominated by $P$.  Thus as elements in $C(X, \Bbb Z)$ we have 
that $\varphi_{*}([P]) - [P]$ is a nonzero function in $-K_{0}^{+}(C(X))$ 
and hence 
$H_{\varphi} \cap K_{0}^{+}(C(X)) \neq 0$.  Thus by Theorem 0.2 we have 
that $C(X)\times_{\varphi}\Bbb Z$ is not AF embeddable.

$(\Leftarrow)$ Assume that $X(\varphi) = X$.  Now take $f \in C(X, \Bbb Z)$ 
and assume that $f \circ \varphi^{-1} - f \geq 0$.  Let 
$\{s_1, \ldots, s_k, -t_1, \ldots, -t_j\}$ be the range of 
$f$, where 
the $s_i, \ t_i$ are all nonnegative integers and $s_i > s_{i + 1}, \ 
t_l > t_{l + 1}$ (and perhaps $s_k = 0$).  Then define the clopen sets 
$E_i = f^{-1}(s_i), \  F_l = f^{-1}(t_l)$ and notice that these sets 
form a pairwise disjoint, finite clopen cover of $X$. Letting 
$P_E$ 
denote the characteristic function of a clopen set, $E$, we can write

\begin{center}
 $f = \sum_{i = 1}^{k}s_iP_{E_i} - \sum_{l = 1}^{j}t_lP_{F_l}$
\end{center}
and thus
\begin{center}
$$f \circ \varphi^{-1} - f = 
  \sum_{i = 1}^{k} s_iP_{\varphi(E_{i})} + \sum_{l = 1}^{j} 
  t_{l}P_{F_l} - \sum_{i = 1}^{k} s_iP_{E_{i}} 
  - \sum_{l = 1}^{j} t_{l}P_{\varphi(F_j)} 
$$
\end{center}

     Now, given $x \in E_1$ we have that the second summation above 
vanishes at $x$.  Also, since $f \circ \varphi^{-1} - f \geq 0$, we have
a unique index $i$  such that 
\begin{center}
 $f \circ \varphi^{-1}(x) - f(x) = s_i - s_1  \geq 0$
\end{center}
However, as $s_1 > s_i$, for $2 \leq i \leq k$, we conclude that 
$s_i = s_1$ and hence $\varphi(E_1) \supset E_1$.  But this implies 
that $\varphi(E_{1}^{c}) \subset E_{1}^{c}$, where $E^{c}$ denotes the 
complement.  But then by hypothesis (and Lemma 4.2) we see that 
$\varphi(E_1) = E_1$.  Repeating this argument we get that 
$\varphi(E_i) = E_i$ for $1 \leq i \leq k$.  An obvious adaptation of this 
argument shows equality for the $F_l$'s.  Hence 
$f \circ \varphi^{-1} - f = 0$, i.e. $H_{\alpha} \cap K_{0}^{+} (C(X)) = 0$. 
$\Box$

\vspace{2mm}

{\bf Remark 4.8} \ \ We may now give a simple example showing that in 
general, (FOP) is not a necessary condition for AF embeddability.

     Let $X$ be the one point compactification of $\Bbb Z$ and define 
$\varphi : X \rightarrow X$ to be the homeomorphism taking 
$n \mapsto n + 1$ and $\infty \mapsto \infty$.  As Pimsner points out in 
[Pi], every point of $X$ will be pseudo-nonwandering (and hence 
$C(X)\times_{\varphi}\Bbb Z$ will be AF embeddable) while it is not 
hard to see that every element of 
$K_0(C(X)) = C(X, \Bbb Z)$ will have infinite orbit under the 
iterates of $\varphi_{*}$.     

\vspace{2mm}

     We now show that the existence of
enough $\alpha_{*}$-invariant states is a sufficient condition for 
AF embeddability.  For convenience we will assume that A is unital and, 
again following the terminology in [Da], we define a 
{\em state} on $K_0(A)$ to be a contractive 
group homomorphism, $\tau$: $K_0(A)$ $\rightarrow  (\Bbb R,\Bbb R_{+}, 
[0, 1])$, with $\tau([1_{A}]) = 1$. 
 There is a 1-1 correspondence between the states on 
$K_0(A)$ and the tracial states on the algebra, $A$ (Theorem IV.5.3 in [Da]).
  Recall 
that a (tracial) state is said to be {\em faithful} if the only positive 
element in it's kernel is the zero element.  We will denote by 
$S_{\alpha}$ the set of all {\em $\alpha_{*}$-invariant} states, i.e. 
those states for which $\tau\circ\alpha_{*} = \tau$.  Note that 
when $A$ is unital, $S_{\alpha}$ is never empty (since $\Bbb Z$ is 
amenable).

     The author would like to thank Professors L.G. Brown and N.C. Phillips 
for pointing out a substantial simplification in the  proof of 
the following theorem.

\vspace{4mm}

{\bf Theorem 4.9} \ \ Let $A$ be unital and  $\alpha \in \Aut (A)$ be given.
Assume that for every projection, $p \in A$, there exists a state, 
$\tau \in S_{\alpha}$, such that $\tau([p]) > 0$. Then 
$A\times_{\alpha}\Bbb Z$ is AF embeddable.  In particular, if $S_{\alpha}$ 
contains a faithful state or 
if $A\times_{\alpha}\Bbb Z$ admits a faithful tracial state then 
it is AF embeddable.

\vspace{2mm}

{\bf Proof} \ \ Again, we prove the contrapositive.  So assume that 
there exists an element $x \in K_0(A)$ and a projection $q \in M_n(A)$ 
such that $\alpha_{*}(x) - x = [q] \neq 0$.  Then clearly 
$\tau([q]) = 0$, for every $\tau \in S_{\alpha}$.  But since 
$\Gamma(A)$ generates $K_{0}^{+}(A)$, we can find projections 
$p_1, \ldots ,p_n \in A$ such that $[q] = [p_1] + \ldots + [p_n]$ and 
thus we see that $\tau([p_i]) = 0$ for every $\tau \in S_{\alpha}$ 
and $1 \leq i \leq n$.  $\Box$

\vspace{2mm}

We are indebted to Kishimoto Sensei for pointing out the following 
corollary.

\vspace{2mm}

{\bf Corollary 4.10} \ \ If $A$ is simple and unital then 
$A\times_{\alpha}\Bbb Z$ is AF embeddable for every $\alpha \in \Aut (A)$.

\vspace{2mm}

{\bf Proof} \ \ This follows from Theorem 4.9 since every state 
on a simple dimension group is faithful.(cf. [Ef], Corollary 4.2) $\Box$ 

\vspace{2mm}

{\bf Remark 4.11} \ \ Corollary 4.10 shows the contrast between 
crossed products
of unital and non-unital algebras.  Indeed, it was first shown in [Cu] 
that the stabilizations of the
Cuntz algebras, $\cal O_{n} \otimes \Bbb K$, are isomorphic to 
crossed products of 
(non-unital) simple AF algebras.  Thus crossed products of non-unital 
simple AF algebras can be purely infinite and hence not AF embeddable.

\vspace{2mm}

     As a final application, we present a few more conditions which 
characterize AF embeddability.

\vspace{4mm}

{\bf Theorem 4.12} \ \ If $A$ is an AF algebra (not necessarily unital) 
and $\alpha \in \Aut (A)$ then
the following are equivalent
\begin{enumerate}
 \item $A\times_{\alpha}\Bbb Z$ is AF embeddable

 \item $\tilde{\rm A}\times_{\tilde{\alpha}}\Bbb Z$ is AF embeddable

 \item For every $m \in \Bbb Z, A\times_{\alpha^{m}}\Bbb Z$ is AF 
       embeddable

 \item There exists an integer, $m \neq$ 0, such that 
       $A\times_{\alpha^{m}}\Bbb Z$ is AF embeddable

 \item There exists an AF algebra, $B$, and a group homomorphism \\
       $\theta: K_{0}( A\times_{\alpha}\Bbb Z) \rightarrow$
       $K_0(B)$ with $i) \  \theta([p]) \in \Gamma(B)$ and $ii) 
        \ [p] \in Ker(\theta) \Rightarrow p = 0$, for every 
       projection $p \in A$.
\end{enumerate}

\vspace{4mm}

{\bf Proof} \ \ $(1)\Leftrightarrow(2)$  It is easy to see that 
\begin{center}
  $K_0(\tilde{A}) = K_0(A)\oplus \Bbb Z$\\
  $H_{\tilde{\alpha}} = H_{\alpha}\oplus 0$
\end{center}

     Thus it is clear that $H_{\tilde{\alpha}} \cap K_{0}^{+}(\tilde{A}) 
= 0 \Leftrightarrow H_{\alpha} \cap K_{0}^{+}(A) = 0$.

$(1) \Leftrightarrow (5)$ The hypotheses of (5) imply that $\theta \circ 
i_{*}$ is a faithful homomorphism (where $i: A \hookrightarrow 
A\times_{\alpha}\Bbb Z$ is the natural 
inclusion)       with $\alpha_{*}(x) - x \in 
Ker(\theta \circ i_{*}$), for all $x \in \Gamma(A)$.  However, since 
$\Gamma(A)$ generates $K_0(A)$ we see that $H_{\alpha} \subset 
Ker(\theta\circ i_{*})$.  Thus the faithfulness of $\theta \circ 
i_{*}$ implies that $H_{\alpha} \cap K_{0}^{+}(A) = \{0\}$.

    The converse is  trivial.

\vspace{2mm}

    By the equivalence of (1) and (2), we may assume for the remainder of 
the proof that A is unital.

\vspace{2mm}

(1)$\Rightarrow$(3) \ \ If $u \in A\times_{\alpha}\Bbb Z$ and $v \in 
A\times_{\alpha^{m}}\Bbb Z$ are the distinguished unitaries, then it is
routine to check that the covariant representation

\begin{center} $v \longmapsto u^{m}$ \end{center}
\begin{center} $a \longmapsto a$       \end{center}
defines an embedding:  $A\times_{\alpha^{m}}\Bbb Z$ $\hookrightarrow$
$A\times_{\alpha}\Bbb Z$.

(3) $\Rightarrow$ (4) is immediate

(4) $\Rightarrow (1)$ \ \ Assume $\varphi: A\times_{\alpha^{m}}\Bbb Z$
$\longrightarrow B$ is an embedding into an AF algebra, $B$.  Notice that in
$K_0(B)$ we have that for every projection, $p \in A, [\varphi(p)] = 
[\varphi(upu^{*})] = [\varphi(\alpha^{m}(p))]$.  We will assume that
m is positive, for if m is negative, it will be clear how to adapt the
following argument.

     Now, define $B_{m} = \oplus_{0}^{m - 1}B$ and let $\beta \in 
\Aut (B_{m}$) be the backwards cyclic permutation of the summands of 
$B_{m}$.  That is, $\beta(b_{0}\oplus \cdots \oplus b_{m - 1}) = 
b_{1}\oplus \cdots \oplus b_{m -1} 
\oplus b_{0}$.  Now, identifying $A$ with it's image in $A\times_{\alpha^{m}}
\Bbb Z$, we define $\psi_{i} = \varphi|_{A}\circ\alpha^{i}, 0 \leq i \leq m-1$.
Then we define $\psi = \oplus_{0}^{m - 1} \psi_{i}: A \rightarrow B_{m}$,
and it is clear that $\psi$ is a *-monomorphism.  It is also clear that 
$\beta$ satisfies (FOP).  Finally, using the fact that in $K_0(B)$ we
know [$\varphi(p)] = [ \varphi(\alpha^{m}(p))]$, one easily checks that
$\beta_{*} \circ \psi_{*} = \psi_{*} \circ \alpha_{*}$.  Thus the 
conclusion follows from Theorem 4.2 and Proposition 3.1. $\Box$

\vspace{2mm}

{\bf Remark 4.13} \ \ It follows from Lemma 2.8 in [Vo] that if 
$\cal U = \otimes_{n \geq 1}M_n$ is 
the Universal UHF algebra, then condition (5) implies that the 
crossed product can be embedded into $B\otimes \cal U$. 

\vspace{2mm}

{\bf Remark 4.14} \ \ The author believes the equivalence of (1) 
and (2) to be a nontrivial 
application of Theorem 0.2 as he has been unable to find a direct 
proof of $(1) \Rightarrow (2)$.

\vspace{2mm}

{\bf Remark 4.15} \ \ The equivalence of (1), (3) and (4) was proved in 
[Pi] when $A$ is assumed to be any (not necessarily AF) unital, 
abelian C$^{*}$-algebra.  

\vspace{7mm}

{\large \bf Section 5: K-Theory}

\vspace{7mm}

     In this section we will show that our previous methods can be used 
to get rationally injective maps on $K_0(A\times_{\alpha}\Bbb Z)$ 
(whenever $A\times_{\alpha}\Bbb Z$ is AF embeddable).  That is, we will 
show that if $A\times_{\alpha}\Bbb Z$ is AF embeddable, then one can 
find a *-monomorphism, $\varphi : A\times_{\alpha}\Bbb Z \rightarrow B$, 
where $B$ is AF and $\varphi_{*} : K_0(A\times_{\alpha}\Bbb Z) 
\rightarrow K_0(B)$ is rationally injective.  The author presented the 
main results of Sections 2, 3 and 4 in the Functional Analysis Seminar 
at Purdue University.  He would like to thank Professor Larry Brown 
for asking the questions that led to this section of the paper.

     The following definition is taken from [Vo].

\vspace{2mm}

{\bf Definition 5.1} \ \ If $\beta \in \Aut(B)$ then $\beta$ is called 
$limit \ periodic$ if there exists an increasing nest, $\{B_n\}$, of 
finite dimensional subalgebras (with $ B = \overline{\bigcup B_n}$) 
and a sequence 
of positive integers, $d_n \in \Bbb N$, such that $\beta(B_n) = B_n$ and 
$\beta^{d_n}|_{B_n} = id_{B_n}$.

\vspace{2mm}

     In Lemma 2.8 of [Vo] it is shown that crossed products of AF algebras 
by limit periodic automorphisms are AF embeddable.  The following lemma 
may be well known to the experts, but the author is not aware of a specific 
reference and thus includes a proof for completeness.

\vspace{2mm}

{\bf Lemma 5.2} \ \ If $B$ is unital and $\beta \in \overline{Inn(B)}$ 
is limit 
periodic (with 
respect to $\{B_n\}$ and with integers $\{d_n\}$) then the 
embedding of 
$B\times_{\beta} \Bbb Z$  into $B\otimes \cal U^{\prime}$ 
(where $\cal U^{\prime}$ is a UHF algebra) 
constructed in Lemma 2.8 of [Vo] induces an injective map on 
$K_0(B\times_{\beta} \Bbb Z)$.

\vspace{2mm}

{\bf Sketch of Proof} \ \ (cf. Lemma 2.8 in [Vo]).
 By the Pimsner-Voiculescu six term exact sequence, we have that 
$i_{*} : K_0(B) \rightarrow K_0(B\times_{\beta} \Bbb Z)$ is an isomorphism,
where $i : B \rightarrow B\times_{\beta} \Bbb Z$ is the natural inclusion 
map.  Letting $\varphi : B\times_{\beta} \Bbb Z \rightarrow 
B\otimes \cal U^{\prime}$ be the embedding constructed in [Vo] 
we have the following commutative diagrams.

\begin{center}
$\xymatrix{
    B \ar[r]^i \ar[dr]^{\varphi} 
    & B\times_{\beta} \Bbb Z \ar[d]^{\varphi} \\
    & B\otimes \cal U^{\prime}     }$
\end{center}

\begin{center}
$\xymatrix{
    K_0(B) \ar[r]^{i_{*}} \ar[dr]^{\varphi_{*}} 
    & K_0(B\times_{\beta} \Bbb Z) \ar[d]^{\varphi_{*}} \\
    & K_0(B\otimes \cal U^{\prime})     }$
\end{center}

Thus, it suffices to show that $\varphi_{*} : K_0(B) \rightarrow 
 K_0(B\otimes \cal U^{\prime})$ is an injective map.  Now, as in the proof of 
Lemma 2.8 in [Vo], we define $U_n = M_{d_n}(\Bbb C)$ and take 
maps $\varphi_{n} 
: B_n \rightarrow B_n \otimes U_n$ with commutativity in the diagram

$$
\begin{CD}
B_n @>j_n>> B_{n + 1} \\
@V\varphi_n VV @V\varphi_{n + 1}VV \\
B_n \otimes U_n @>>> B_{n + 1} \otimes U_{n + 1}
\end{CD}
$$
where we don't need to know what map the lower arrow is given by, $j_n$ is 
the natural inclusion map and 
$\varphi_n$ takes $b \mapsto \sum_{j = 1}^{d_n} \beta^{j}(b)\otimes 
e_{j,j}$, where $e_{i,j}$ are the standard matrix units of $M_{d_n}$.
Since $\varphi$ is defined as the limit of the maps $\varphi_n$, we have 
the following commutative diagram

$$
\begin{CD}
K_0(B_n) @>j_{n*}>> K_0(B) \\
@V\varphi_{n*} VV @V\varphi_{*}VV \\
K_0(B_n \otimes U_n) @>\Psi_{n*}>> K_0(B \otimes \cal U^{\prime})
\end{CD}
$$
where $\Psi_n : B_n \otimes U_n \rightarrow B \otimes \cal U^{\prime}$ is the 
induced map and actually, $\cal U^{\prime} = \otimes_{n \geq 1} M_{d_n}$.  
From this diagram we see that it suffices to show the 
following assertion.

\vspace{1mm}

     $Claim$: If $x \in K_0(B_n)$ and $\varphi_n(x) = 0$ then  
$j_{n*}(x) = 0$.

\vspace{1mm}

     So, assume we have $x = [p] - [q] \in K_0(B_n)$ (where $p, \ q$ are 
projections in the matrices over $B_n$) such that 
$\varphi_{n*}(x) = 0$.  Then with the identification 
$K_0(B_n \otimes U_n) = K_0(M_{d_n}(B_n)) = K_0(B_n)$, we have that 
(in $K_0(B_n)$) $\varphi_n(x) =  
\sum_{j = 1}^{d_n} [\beta^{j}(p)] - \sum_{j = 1}^{d_n} [\beta^{j}(q)] = 0$. 
However, by assumption we have that in $K_0(B)$, $[r] = [\beta^{j}(r)]$, for 
every projection, $r$, in the matrices over $B$ and for every 
$j \in \Bbb Z$.  Thus, 
passing to $K_0(B)$ the above formula becomes $d_n([p] - [q]) = 0$ or 
$d_n(j_{n*}(x)) = 0$ and 
thus (since $K_0(B)$ is torsion free) $j_{n*}(x) = 0$.  $\Box$

\vspace{4mm}

     We now present an analogue of the lemma of Spielberg (Lemma 3.4) which 
gives more control on K-theory.  The main idea in the proof is similar 
to the proof of Lemma 1.14 is [Sp].

\vspace{4mm}

{\bf Lemma 5.3} \ \ If $(G, \ G^{+})$ is a torsion free, ordered group then 
there  exists a cone $\tilde{G}^{+} \supset G^{+}$ such that 
$(G, \ \tilde{G}^{+})$ is totally ordered and hence is a dimension 
group.  

\vspace{2mm}

{\bf Proof} \ \ Let $\cal L = \{ E \subset G : i) \ E \supset G^{+} \ \ 
ii) \  E \text{ is a semigroup} \ \ iii) \ E\cap (-E) = 0 \ \}$ be 
partially ordered by inclusion.  It is easy to check that the hypotheses 
of Zorn's Lemma are satisfied and thus we may choose a maximal element, 
$\tilde{G}^{+}  \in \cal L$.  Now, choose $x \in G\backslash \tilde{G}^{+}$ 
and consider the 
semigroup $\{ nx + e : n \in \Bbb N, \ e \in \tilde{G}^{+} \}$.  
This semigroup 
clearly (properly) contains $\tilde{G}^{+} $, and thus by maximality 
we get that 
there exist 
nonnegative integers, $n, \ m \in \Bbb N$ (at least one of which is nonzero)
, and elements, $e, \ f \in \tilde{G}^{+} $ , 
such that $nx + e = -(mx + f) \neq 0$.  Thus $(n + m)x \in -\tilde{G}^{+} $.

     We now claim that $(G, \ \tilde{G}^{+})$ is unperforated.  
So assume that we 
have found some nonzero element, $x \in G$, such that $kx \in \tilde{G}^{+}$ 
while $x \notin \tilde{G}^{+} $.  Then the argument above 
(since $x \notin \tilde{G}^{+} $) provides us 
with a positive integer (letting $j = n + m$, from above) such that 
$jx \in -\tilde{G}^{+} $.  Thus, $(kj)x \in \tilde{G}^{+} \cap -\tilde{G}^{+}
 = 0$.  
But this is a 
contradiction to the hypothesis that $G$ is torsion free and hence 
$(G, \ \tilde{G}^{+})$ is unperforated. 

     However, if we now take $x \in G\backslash \tilde{G}^{+} $ then 
we can find 
a positive integer $j$ such that $jx \in -\tilde{G}^{+} $ which implies 
(since $(G, \ \tilde{G}^{+} )$ is unperforated) that $x \in -\tilde{G}^{+} $.
  Thus $(G, \ \tilde{G}^{+} )$ is 
totally ordered.  $\Box$

\vspace{4mm}

{\bf Remark 5.4} \ \ If $A\times_{\alpha}\Bbb Z$ is AF embeddable then it 
is easy to see that $K_0(A\times_{\alpha}\Bbb Z) = K_0(A)/H_{\alpha}$ is 
an ordered group with the cone $K_{0}^{+}(A) + H_{\alpha}$.  Clearly, the 
only thing that needs to be checked is that $K_{0}^{+}(A) + H_{\alpha} \cap 
-(K_{0}^{+}(A) + H_{\alpha}) = \{0\}$.  However, this follows easily from 
the fact that $K_{0}^{+}(A) \cap H_{\alpha} = \{0\}$ by assumption (and 
Theorem 0.2).

\newpage

{\bf Theorem 5.5} \ \  Assume that $A\times_{\alpha}\Bbb Z$ is AF embeddable 
and $H$ is a subgroup of $K_0(A\times_{\alpha}\Bbb Z) = K_0(A)/H_{\alpha}$ 
with the property 
that $H \cap (K_{0}^{+}(A) + H_{\alpha}) = \{0\}$ and 
$K_0(A\times_{\alpha}\Bbb Z)/H$ is torsion free.  Then one may choose 
the AF embedding of $A\times_{\alpha}\Bbb Z$ such that the kernel of 
the induced map on $K_0(A\times_{\alpha}\Bbb Z)$ is precisely $H$.

\vspace{2mm}

{\bf Proof} \ \ To ease our notation somewhat, we begin by defining 
\begin{center}
 $G = K_0(A\times_{\alpha}\Bbb Z)/H$\\
 $G^{+} = (K_{0}^{+}(A) + H_{\alpha}) + H$
\end{center}
where $K_{0}^{+}(A) + H_{\alpha}$ was shown to be a cone of 
$K_0(A\times_{\alpha}\Bbb Z) = K_0(A)/H_{\alpha}$ in Remark 5.4 and thus 
$G^{+}$ is just the natural image of this cone in $G$.  Now we claim that 
$(G, \ G^{+})$ is an ordered group.  It is clear that we have $G^{+} + 
G^{+} \subset G^{+}$ and $G = G^{+} - G^{+}$ and thus we only have to 
check that $G^{+} \cap -(G^{+}) = \{0\}$.  So assume there are elements 
$x, \ y \in K_{0}^{+}(A)$ such that $(x + H_{\alpha}) + H = -(y + H_{\alpha}) 
+ H$ (in G).  This implies that $(x + y) + H_{\alpha} \in H \cap 
(K_{0}^{+}(A) + H_{\alpha}) = \{0\}$ (in $K_0(A\times_{\alpha}\Bbb Z)$)  
and hence $x + y \in H_{\alpha} \cap 
K_{0}^{+}(A) = \{0\}$ (in $K_0(A)$).  But, this implies that 
$x = y = 0$.

     Thus we have that $(G, \ G^{+})$ is a torsion free, ordered group and 
hence by Lemma 5.3 we can find a cone $\tilde{G}^{+} \supset G^{+}$ and 
an AF algebra, $B$, such that $(K_0(B), \ K_{0}^{+}(B), \ \Gamma(B)) = 
(G, \ \tilde{G}^{+}, \ \tilde{G}^{+})$.  Note that if we let 
$\pi : K_0(A) \rightarrow K_0(B) = G$ be the canonical projection 
then $\pi$ is a contractive, faithful group homomorphism.  Hence, by 
Fact 1.1, there exists a *-monomorphism, $\varphi : A \rightarrow B$ with 
$\varphi_{*} = \pi$.  The crucial properties in the remainder of the proof 
will be, 

\begin{enumerate}
 \item $\varphi : A \rightarrow B$ is a unital *-monomorphism
 \item $K_0(B) =  K_0(A\times_{\alpha}\Bbb Z)/H$
 \item $Ker(\varphi_{*}) = Ker(\pi)$ where $\pi$ is the composition 
       of the canonical projection maps $K_0(A) \rightarrow 
       K_0(A\times_{\alpha}\Bbb Z) \rightarrow 
       K_0(A\times_{\alpha}\Bbb Z)/H$.

\end{enumerate}

     We now show that we may assume that $A$ and $B$ are unital and 
that $\varphi : A \rightarrow B$ is a unital injection.  By the split 
exactness of the sequence 

\begin{center}
$0 \rightarrow A\times_{\alpha}\Bbb Z \rightarrow 
\tilde{A}\times_{\tilde{\alpha}}\Bbb Z \rightarrow C(\Bbb T) 
\rightarrow 0$
\end{center}
we have that  
$K_0(\tilde{A}\times_{\tilde{\alpha}}\Bbb Z) = 
K_0(A\times_{\alpha}\Bbb Z) \oplus \Bbb Z$ and the image of $H$ under 
this identification is $\tilde{H} = H \oplus 0  
\subset K_0(\tilde{A}\times_{\tilde{\alpha}}\Bbb Z)$.  We now claim that 
properties 1), 2) and 3) above still hold with the unitizations of $A, \ 
B$ and $\varphi$.

     Now we recall that $K_0(\tilde{A}) = K_0(A) \oplus \Bbb Z$ and 
$K_0(\tilde{B}) = K_0(B) \oplus \Bbb Z$ and hence $ker(\tilde{\varphi}_{*}) 
= ker(\varphi_{*}) \oplus 0 \subset K_0(A) \oplus \Bbb Z$.  
Also, we observe that 

\begin{center}
$K_0(\tilde{A}\times_{\tilde{\alpha}}\Bbb Z)/\tilde{H} = 
(K_0(A\times_{\alpha}\Bbb Z) \oplus \Bbb Z)/(H \oplus 0) = 
K_0(B) \oplus \Bbb Z = K_0(\tilde{B})$
\end{center}
with all these identifications being natural.  Thus properties 1) and 
2) are satisfied. Finally, it is easy to see that the kernel of the 
composition of the projections maps 

\begin{center}
$K_0(\tilde{A}) \rightarrow K_0(\tilde{A}\times_{\tilde{\alpha}}\Bbb Z) 
\rightarrow K_0(\tilde{A}\times_{\tilde{\alpha}}\Bbb Z)/\tilde{H}$ 
\end{center}
is precisely $Ker(\varphi_{*}) \oplus 0 = Ker(\tilde{\varphi}_{*})$ since 
the above sequence is really just 
\begin{center}
$K_0(A) \oplus 0 \rightarrow K_0(A\times_{\alpha}\Bbb Z) \oplus \Bbb Z 
\rightarrow (K_0(A\times_{\alpha}\Bbb Z) \oplus \Bbb Z)/(H \oplus 0)$
\end{center}

     Thus it suffices to prove the theorem with 
$\tilde{A}\times_{\tilde{\alpha}}\Bbb Z$ and $\tilde{H}$ and hence we may 
assume that $A$ and $B$ are unital and $\varphi$ is a unit preserving, 
injective *-homomorphism with properties 1), 2) and 3) above.

     Now, by the proof 
of Proposition 3.1 (and Proposition 2.8) we may further  assume 
that, $\varphi : A \rightarrow B\otimes \cal U$, where $\cal U$ is the 
Universal UHF 
algebra and we have commutativity in the diagram 

$$
\begin{CD}
A @>\varphi>>B\otimes \cal U \\
@VAdv\circ\alpha VV  @VVAdu\circ(id_{B}\otimes\sigma) V \\
A @>\varphi>>B\otimes \cal U
\end{CD}
$$
where $v \in A$, $u \in B\otimes \cal U$ are unitaries and $\sigma \in 
\Aut(\cal U)$ is the automorphism (with the Rohlin property) from Example 2.2.
Note that it was necessary to first arrange that $\varphi$ be a unital 
map in order to appeal to Proposition 2.8 and that we have not 
changed the kernel of $\varphi_{*}$ in doing so (although  we now have 
that $K_0(B) = K_0(A\times_{\alpha}\Bbb Z)/H$ sits injectively inside
 $K_0(B\otimes \cal U)$).
Note also that $id_{B}\otimes\sigma$ is limit periodic.

Now, from the Pimsner-Voiculescu six term exact sequence we have that 
$K_0(B\otimes\cal U) = K_0(B\otimes\cal U 
\times_{id_{B}\otimes\sigma}\Bbb Z)$.  Thus commutativity in the diagram 

$$
\begin{CD}
A @>i>> A\times_{Adv\circ\alpha}\Bbb Z \\
@V\varphi VV @VV\tilde{\varphi} V \\
B\otimes\cal U @>>> B\otimes\cal U \times_{Adu\circ(id_{B}\otimes\sigma)}\Bbb Z
\end{CD}
$$
implies commutativity in the diagram
$$
\begin{CD}
K_0(A) @>i_{*}>> K_0(A\times_{Adv\circ\alpha}\Bbb Z) \\
@V\varphi_{*} VV @VV\tilde{\varphi}_{*} V \\
K_0(B\otimes\cal U) @>\cong>> K_0(B\otimes\cal U 
\times_{Adu\circ (id_{B}\otimes\sigma)}\Bbb Z)
\end{CD}
$$
Where the $\tilde{\varphi}$ on the right side of the first diagram is 
now the natural extension of $\varphi$ to the crossed products (i.e. it 
no longer denotes the unital extension) and $i$ is 
the natural inclusion map.

     Finally, since  $Ker(\varphi_{*}) = 
Ker(\tilde{\varphi}_{*}\circ i_{*}) = Ker(\pi)$ where $\pi$ was the 
composition of the maps 
 
\begin{center}
$K_0(A) \rightarrow 
       K_0(A\times_{\alpha}\Bbb Z) \rightarrow 
       K_0(A\times_{\alpha}\Bbb Z)/H = K_0(B) \hookrightarrow 
K_0(B\otimes\cal U)$
\end{center}       
it is now a routine exercise to 
verify that $Ker(\tilde{\varphi}_{*}) = H$ (under the identifications 
$K_0(A\times_{Adv\circ\alpha}\Bbb Z) = K_0(A\times_{\alpha}\Bbb Z)$ and 
$K_0(B\otimes\cal U \times_{Adu(\circ id_{B}\otimes\sigma)}\Bbb Z)  = 
K_0(B\otimes\cal U \times_{id_{B}\otimes\sigma}\Bbb Z)$).  Thus the 
proof is complete since $id_{B}\otimes\sigma$ is limit periodic and 
hence from Lemma 5.2 we have that 
$K_0(B\otimes\cal U \times_{id_{B}\otimes\sigma}\Bbb Z)$ gets embedded into 
$K_0(B\otimes\cal U\otimes \cal U^{\prime})$ where $\cal U^{\prime}$ is 
the UHF algebra  
constructed in Lemma 5.2. $\Box$ 

\vspace{2mm}

{\bf Corollary 5.6} \ \ If $A\times_{\alpha}\Bbb Z$ is AF embeddable, then 
one may choose an embedding which  induces a rationally injective map on 
$K_0(A\times_{\alpha}\Bbb Z)$.

\vspace{2mm}

{\bf Proof} \ \ Let $T(A\times_{\alpha}\Bbb Z) \subset 
K_0(A\times_{\alpha}\Bbb Z) = K_0(A)/H_{\alpha}$ denote the torsion subgroup.
  To apply 
Theorem 5.5 we only need to see that $T(A\times_{\alpha}\Bbb Z) 
\cap (K_{0}^{+}(A) + H_{\alpha}) = \{0\}$.  But this follows easily 
from the fact that $K_{0}^{+}(A) \cap H_{\alpha} = \{0\}$.  $\Box$ 
  
\vspace{2mm}

{\bf Remark 5.7} \ \ It is easy to prove that if $K_0(A)$ is a divisible 
group then $K_0(A\times_{\alpha}\Bbb Z)$ is always torsion free.  However, 
it not hard to construct an example where $A\times_{\alpha}\Bbb Z$ is 
AF embeddable and $K_0(A\times_{\alpha}\Bbb Z)$ has torsion.

\vspace{2mm}

{\bf Remark 5.8} \ \ It is clear that, in general, our constructions will 
not yield isomorphisms on K-theory since we must tensor with some UHF 
algebra to get the desired embeddings.

\end{document}